\documentclass{amsart}

\usepackage{amssymb}

\theoremstyle{plain}

\numberwithin{equation}{section}

\newtheorem*{conj}{Conjecture}

\newtheorem*{teo}{Theorem A}
\newtheorem*{teob}{Theorem B}

 \newtheorem*{cor}{Corollary C}

\newtheorem{hypo}[equation]{Hypotheses}

\newtheorem{lem}[equation]{Lemma}

\newcommand{\Irr}{\operatorname{Irr}}

\newcommand{\GF}{\operatorname{GF}}
\newcommand{\dl}{\operatorname{dl}}

\newcommand{\core}{\operatorname{core}}
\newcommand{\Ker}{\operatorname{Ker}}

\theoremstyle{definition}

\newtheorem{example}[equation]{Example}

\begin{document}	
\title{Derived Length and Products of Conjugacy Classes}

\author{Edith Adan-Bante}

\address{University of Southern Mississippi Gulf Coast, 730 East Beach Boulevard,
 Long Beach MS 39560}

\email{Edith.Bante@usm.edu}



\date{2006}

\begin{abstract}
Let $G$ be a supersolvable group and $A$ be a conjugacy class
of $G$. Observe that for some integer 
$\eta(AA^{-1})>0$, $AA^{-1}=\{a b^{-1}\mid a,b\in A\}$ is the union of $\eta(AA^{-1})$
distinct conjugacy classes of $G$. Set ${\bf C}_G(A)=\{g\in G\mid a^g=a\text{ for all } a\in A\}$. Then the derived length of $G/{\bf C}_G(A)$ is less or
equal than $2\eta(A A^{-1})-1$.
\end{abstract}
\maketitle
\begin{section}{Introduction}

 Let $G$ be a finite group, $A$ be 
a conjugacy class of $G$ and $e$ be the identity of $G$.
Let $X$ be a nonempty $G$-invariant subset of $G$, i.e. 
$X^g=\{ g^{-1} xg\mid x\in X\}=X$ for all $g\in G$.
  Then for some integer $n>0$, $X$ is  a union of 
  $n$ distinct conjugacy classes of $G$. Set
 $\eta(X)=n$.
 
  We can check that given any two conjugacy classes
 $A$ and $B$ of $\,G$, the product $A B=\{ab \mid a\in A  , b\in B\}$ of
 $A$ and $B$ is a 
 $G$-invariant set. Thus $\eta(A B)$ is the number of distinct conjugacy 
 classes of $G$ such that $A B$ is the union of those classes.
 
Denote by ${\bf C}_G(A)=\{g\in G \mid a^g =a \text{ for all } a\in A\}$ the centralizer of the set $A$ in $\,G$. If $\,G$ is a
solvable group, denote by $\dl(G)$ the derived length of $G$. Let 
${\bf Z}(G)$ be the center of $\,G$.

In this note, we are exploring the relations between the structure of the group $G$ and 
the product $A B$ of  some conjugacy 
 classes $A$ and $B$ of $\,G$. More specifically,
 we are exploring the relation  between the derived length of some section of $\,G$ and
 properties of $A B$.
 
\begin{teo}
Let $G$ be a finite group, $A$ and $B$ be conjugacy classes of $\,G$ and $a\in A$, $b\in B$.
  If  $a, b\in {\bf Z}_2(G)$, i.e. $[a,g]\in {\bf Z}(G)$ for all
  $g\in G$, then 
${\bf C}_G(A)\cap {\bf C}_G(B)\supseteq [G,G]$. In particular $\dl(G/{\bf C}_G(A))\leq 1$.
\end{teo}

Given a finite solvable group $G$ and conjugacy classes $A$ and $B$ of $\,G$,
is there any relationship between the derived length of $G$ and 
 $\eta(A B)$? In general, the answer seems to be no. For instance
$A \{e\}=A$ for any finite group $G$ and any conjugacy class $A$ of $G$. 
 Thus $\eta(A B)$ may not give us information
 about $\dl(G)$, but it does give us a linear bound on $\dl(G/
 {\bf C}_G(A))$
 when $B=A^{-1}$ and $G$ is supersolvable. More precisely
 
 \begin{teob} For any finite supersolvable group $G$ and any conjugacy class 
 $A$ of $\,G$ we have that 
 \begin{equation}
 \dl(G/{\bf C}_G(A)) \leq 2 \eta(A A^{-1}) -1.
 \end{equation} 
 \end{teob}
 
 An application of this result is the following
 
 \begin{cor}
 For any finite supersolvable group and any conjugacy classes 
 $A$, $B$ of $\,G$ such that $A B\cap {\bf Z}(G)\neq \emptyset$,
  we have that 
 \begin{equation*}
 \dl(G/{\bf C}_G(A)) \leq 2 \eta(A B) -1.
 \end{equation*} 
 \end{cor}
 
 We want now to point out the ``dual" situation with characters. 
 
 Denote by $\Irr(G)$ the set of irreducible complex characters of $G$. We can check that
  the product of characters is a character. Therefore $\chi\psi$ is a character of $G$ for any 
 $\chi,\psi\in \Irr(G)$. It is known that a character can be expressed as an integral linear 
 combination of irreducible characters. Then the decomposition of the character $\chi\psi$
 into its distinct irreducible constituents $\theta_1, \theta_2, \ldots, \theta_n$ has the 
 form
 \begin{equation*}
 \chi\psi=\sum_{i=1}^n m_i \theta_i
 \end{equation*} 
 \noindent where $n>0$ and $m_i$ is the multiplicity of $\theta_i$. 
 Set $\eta(\chi\psi)=n$, so that $\eta(\chi\psi)$ is the number of distinct irreducible 
 constituents of the product $\chi\psi$. 
Define $\overline{\chi}(g)$  to be the complex conjugate $\overline{\chi(g)}$
 of $\chi(g)$ for all $g\in G$. 
 
 In Theorem A of \cite{length1}, it is proved that there exist constants 
 $c$ and $d$ such that for any finite 
 solvable group $G$ and any irreducible character $\chi$ of $G$, 
 \begin{equation*}
 \dl(G/\Ker(\chi)) \leq c\eta(\chi\overline{\chi}) +d.
 \end{equation*}
 \noindent If, in addition, $G$ is a supersolvable group, then we may take $c=2$ and $d=-1$.
 
 We regard ${\bf C}_G(A)$ for the conjugacy class $A$ as the dual 
 of $\Ker(\chi)$ for a character $\chi\in \Irr(G)$, the conjugacy class $A^{-1}$ for 
 $A$ as the dual of  $\overline{\chi}$ for the character $\chi\in \Irr(G)$.
   Thus we regard Theorem B as the 
 dual in conjugacy classes of Theorem A  of \cite{length1} for supersolvable groups. 
 In light of Theorem A of \cite{length1}, we wonder 

 \begin{conj}
There exist universal constants $q$ and $r$ such that for any 
finite solvable group and any conjugacy class $A$ of $\,G$, we have that
\begin{equation*}
 \dl(G/{\bf C}_G(A)) \leq q \eta(A A^{-1}) + r.
 \end{equation*} 
 \end{conj}
 We want to remark that there are several results showing the ``duality" of 
 products of conjugacy classes and 
 products of characters. For example,  see \cite{arad}, 
 \cite{productchar2} and \cite{productc}, \cite{squares} and \cite{producth},
 \cite{productchar1} and \cite{productcp2}.
 However, we also want to remark that there are results in
 products of conjugacy classes that do not hold true for 
 the ``dual"
 in  characters.
  For instance, it has been proved that 
  the product of nontrivial  conjugacy classes
 in $A_n$, for $n\geq 5$, is never a conjugacy classes [see \cite{aradbook}]. On the other hand,
  if $n\geq 5$ is a perfect square, there exist irreducible characters $\chi$ and $\psi$ in $A_n$ such that the product $\chi\psi$  is also an irreducible
  [see \cite{zisser}].

{\bf Acknowledgment.} I would like to thank Professor Everett C. Dade for his
suggestions to improve both the result and the presentation of this note.
 I would like also to thank FEMA for providing me with temporary housing in the 
 aftermath of 
  hurricane Katrina. 
\end{section}
 
 \begin{section}{The function $\eta$}
 Let $G$ be a finite group. Observe that if 
 $A$ is a conjugacy class of $G$, then for any $a\in A$
 we have that $A=a^G=\{a^g\mid g\in G\}$.
 
 In this section, we show that given a subgroup $H$ with ${\bf C}_G (a)\subseteq H\subseteq G$, we may 
 not have a relation between 
 $\eta(a^G (a^{-1})^G)$ and $\eta(a^H (a^{-1})^H)$.
 
\begin{example} 
If $H$ is a subgroup of $G$ with  ${\bf C}_G (a)\subset H\subseteq G$,  
 then we need not have     
$\eta(a^G (a^{-1})^G)\geq \eta(a^H (a^{-1})^H)$. 
\end{example}
\begin{proof}
	Fix a prime $p$. 
Let $\GF(p)=\{0,1,\ldots,p-1\}$ be a finite field with $p$ elements. 
Denote by $F^{*}$ the group of units of  $\GF(p)$.  
Also denote by $F$ the additive group of  $\GF(p)$.
 Observe that $F^{*}$ acts on $F$ by multiplication. 
	Define $M$ to be the semi-direct product 
$F \rtimes F^{*}$ 
of $F$ by  $F^{*}$. 

Let $C=<c>$ be a cyclic group of order $p$ and  $1_C$ be the identity of $C$. 
Let $K$ be the direct product of $p$-copies of
$C$, i.e  
$$K= C \times \cdots \times C.$$
Thus $K$ is an elementary abelian group of order 
$p^p$. We can check 
 that $$(fm)(x)= f(mx),$$ for any $m \in M$ and $x \in \GF(p)$, defines
an action of $M$ on $K$. 
	 
	Let $G$ be the wreath product of $C$ and $M$ relative to $\GF(p)$,
i.e. $G = K \rtimes M$ and $H= K\rtimes \GF(p)$. 

Set $a=(c,1_C,\ldots, ,1_C)$ in $K$. By Theorem A of \cite{productc}, we have that
$\eta(a^H (a^{-1})^H)=p$. By Proposition 5.4 of  \cite{productc}, we have  that
$\eta(a^G (a^{-1})^G)=2$.
\end{proof}
\begin{example}
If $H$ is a subgroup of $G$ with  ${\bf C}_G (a)\subset H\subseteq G$,  
 then we need not have     
$\eta(a^G (a^{-1})^G)\leq \eta(a^H (a^{-1})^H)$. 
\end{example}
\begin{proof}
Let $G$ be an extra-special group of exponent $p$ and order $p^3$, for some
odd prime $p$. Let $a\in G\setminus {\bf Z}(G)$ and $H={\bf C}_G(a)$. We can check
that $\eta(a^H (a^{-1})^H)=1$ and $\eta(a^G (a^{-1})^G)=p$. 
\end{proof}
\end{section}

\begin{section}{Proof of Theorem A}
{\bf Notation.} Let $G$ be a group and $a\in G$.  Set $[a,G]=\{[a,g]\mid g\in G\}$. 
Observe that $a^G=\{a[a,g]\mid g\in G\}=a[a,G]$.

The following is a minor modification of Lemma 2.1 of \cite{producth}.
\begin{lem}\label{lema1}
  Let $G$ be a finite group,  $N$ be a normal subgroup of $G$ and $a,b \in G$. Then $a^N b^N =ab[a,N]^b[b,N]$.
\end{lem}
\begin{lem}\label{theoa}
Let $G$ be a finite group, $N$ be a normal subgroup of $\,G$ and
$A$ and $B$ be conjugacy classes of $G$. Suppose that 
there exist $a\in A$ and $b\in B$ such that  
  $AB\subseteq ab {\bf Z}(N)$.  Then 
${\bf C}_N (A)\cap {\bf C}_N(B) \supseteq [N,N]$. In particular $\dl(N/{\bf C}_N(A))\leq 1$.
\end{lem}
\begin{proof}
By Lemma \ref{lema1}, we have that $[a,N]^b [b,N]\subseteq {\bf Z}(N)$. Since
$e\in [b,N]$ and $e\in [a,N]$, we have that $[a,N]^b \subseteq {\bf Z}(N)$ and
$[b,N]\subseteq {\bf Z}(N)$. Since $N$ is normal in $G$,
 so is 
${\bf Z}(N)$. Thus  
$[a,N]\subseteq ({\bf Z}(N))^{b^{-1}}={\bf Z}(N)$.

Since $[a,N]\subseteq {\bf Z}(N)$, the map $n\mapsto [a,n]$ is an homomorphism of the group
$N$ into ${\bf Z}(N)$. This homomorphism is trivial on $[N,N]$. So $[N,N]$ centralizes $a$.
Similarly, $[N,N]$ centralizes $b$ and the result is proved.
\end{proof}
\begin{proof}[Proof of Theorem A]
Observe that if $a,b\in {\bf Z}_2(g)$, then
 $AB\subseteq ab {\bf Z}(G)$. The result then follows by
 the previous Lemma.
 \end{proof} 
\end{section}
\begin{section}{Proof of Theorem B}

\begin{hypo}\label{hypo} Let $G$ be a solvable group, $N\trianglelefteq G$ be a normal 
subgroup
of $G$, $A$ be a conjugacy class of $G$. Fix $a\in A$.
 Set $C={\bf C}_G(a)$. Since $a\in A$, observe that 
 $a^G=A$ and $C_G(A)=\core_G(C)$. Set
 $\bar{K}= (KN)/N$ for any subgroup of $K$ of  $\,G$, and 
 $\bar{g}=gN$ for any element $g\in G$. Set $C_N=\{g\in G\mid [a,g]\in N\}$. 
\end{hypo}

  \begin{lem}\label{lema2} 
Assume Hypothesis \ref{hypo}.
 Then the set $C_N$ is a subgroup of $G$ containing $CN$, $C_N/N={\bf C}_{\bar{G}}(\bar{a})$, and  
  \begin{equation}\label{equationH}
  \eta(\bar{a}^{\bar{G}} ({\bar{a}}^{-1})^{\bar{G}})+ \eta((a^{C_N} (a^{-1})^{C_N})^G)-1 \leq \eta(a^G (a^{-1})^G).
  \end{equation}
  \end{lem}
  \begin{proof}
  Let $\rho: G\rightarrow G/N$ be the homomorphism defined
  by $\rho(g)=\bar{g}$.
  We can check that
   $\rho(g)=\bar{g}\in {\bf C}_{\bar{G}}(\bar{a})$ if and only if $[a,g]\in N$, that is if
   and only if $g\in C_N$, i.e. $\rho( a^g a^{-1})=\bar{e}$ if and 
   only if 
   $g\in C_N$. Thus $C_N$ is the inverse image of the group 
   ${\bf C}_{\bar{G}}(\bar{a})$ under $\rho$.
   
   Since
   $\rho( a^g a^{-1})=\bar{e}$ if and only if 
   $g\in C_N$, we have 
  $$( \cup_{g\in C_N} (a^g a^{-1})^G)\cap (\cup_{g\in G\setminus C_N}(a^g a^{-1})^G)=\emptyset,$$ 
  \noindent and $$ \eta(\bar{a}^{\bar{G}} ({\bar{a}}^{-1})^{\bar{G}})-1\leq   \eta( \cup_{g\in G\setminus C_N} (a^g a^{-1})^G).$$
  Observe that 
  $$a^G (a^{-1})^G=  (\cup_{g\in  C_N} (a^g a^{-1})^G) \cup 
  (\cup_{g\in G \setminus C} (a^g a^{-1})^G).$$
  \noindent Since $\cup_{g\in  C_N} (a^g a^{-1})^G=(a^{C_N} (a^{-1})^{C_N})^G$,
  \eqref{equationH} follows.
   \end{proof}
   
   \begin{lem}\label{intersection}
    Let $H$ and $K$ be subgroups of $G$ with $K\subseteq H$. Then 
   $$\dl(\core_G(H)/\core_G(K))\leq \dl(H/\core_H(K)).$$
   \end{lem}
   \begin{proof}
  Since
   $\core_G(H)=\cap_{g\in G} H^g$, the map
   $\varphi: \core_G(H)\rightarrow \oplus_{g\in G} H^g/\core_{H^g}(K^g)$ defined by 
   $\varphi(h)=\oplus_{g\in G} h \core_{H^g}(K^g)$
   is well defined. Observe that the kernel of $\varphi$
   is 
   $\{h\in \core_G(H) \mid h \in K^g \text{ for all }g\in G\}=
\core_G(K)$. Thus  $\core_G(H)/\core_G(K)$ is isomorphic to 
a section of $\oplus_{g\in G} H^g/\core_{H^g}(K^g)$ and so 
$\dl(\core_G(H)/\core_G(K))\leq \dl( \oplus_{g\in G} H^g/\core_{H^g}(K^g))$.

 Since $H^g/\core_{H^g}(K^g)$ is isomorphic to $H/\core_H(K)$,
 then 
 $$ \dl(H/\core_H(K))=\dl(\oplus_{g\in G} H^g/\core_{H^g}(K^g)),$$
 \noindent and the result follows.
   \end{proof}
   
   \begin{lem}\label{lema3} Assume Hypothesis \ref{hypo}. Assume
   also that $N$ is abelian.
  Then 
 \begin{equation}\label{equationequal}
 \dl(C_N/\core_{C_N}(C))\leq \dl(C_N/(C_N\cap {\bf C}_G(N)))+1.
 \end{equation}
 
 If, in addition, we have that $N$ is a cyclic subgroup, then 
 \begin{equation}\label{supersolvable}
 \dl(C_N/\core_{C_N}(C))\leq 2.
 \end{equation}
 
 \end{lem}
 \begin{proof}
 Because $C_N=\{g\in G\mid [a,g]\in N\}$, we have  
 \begin{equation}\label{contain}
 a^{C_N} (a^{-1})^{C_N}\subseteq N.
 \end{equation}
 Set $H= C_N\cap {\bf C}_G(N)$.
 Observe that ${\bf Z}(H)={\bf Z}(C_N\cap {\bf C}_G(N))\supseteq N$. By \eqref{contain} and Lemma \ref{theoa},
  we have that ${\bf C}_H(a)\supseteq [H,H]$. Since
  $H\trianglelefteq C_N$,  $[H,H]\trianglelefteq C_N$. Because 
  $[H,H]\trianglelefteq C_N$ and ${\bf C}_H(a)=C\cap H\supseteq [H,H]$, we have that
   $\core_{C_N}(C)\cap {\bf C}_G(N)=\core_{C_N} (C)\cap H \supseteq [H,H]$.
   Thus 
  \begin{equation}\label{equation2}
 \dl(H/  (H\cap \core_{C_N}(C)))    \leq 1.
 \end{equation}
 Then 
 \begin{equation}\label{counting}
\begin{split}
 \dl(C_N/\core_{C_N}(C))  &\leq \dl(C_N/H) + \dl(H/  (\core_{C_N}(C)\cap H))\\
  &\leq \dl(C_N/H)+1,
 \end{split}
 \end{equation}
 \noindent where the last inequality 
  follows from \eqref{equation2}.
 
 If $N$ is a cyclic group, then the group $G/{\bf C}_G(N)$ is abelian and so 
 $C_N/H=C_N/(C_N\cap {\bf C}_G(N))$ is abelian. Thus \eqref{supersolvable} 
 follows from \eqref{equationequal}
 \end{proof}
  \begin{lem}\label{equal} 
 Let $N$ be a chief factor of $G$ and $a\in G$. Suppose that 
 $$\eta(\bar{a}^{\bar{G}} ({\bar{a}}^{-1})^{\bar{G}})= \eta(a^G (a^{-1})^G).$$
 \noindent  Then 
 $N\subseteq {\bf C}_G(a)=C_N$. In particular
 $$\dl(G/\core_G({\bf C}_G(a)))=\dl(\bar{G}/\core_{\bar{G}}({\bf C}_{\bar{G}}(\bar{a}))).$$
 \end{lem}
 \begin{proof}
 Since $\eta(\bar{a}^{\bar{G}} ({\bar{a}}^{-1})^{\bar{G}})= \eta(a^G (a^{-1})^G)$, 
 by Lemma \ref{lema2}
 we have that $$\eta((a^{C_N} (a^{-1})^{C_N})^G)=1.$$ Since $e\in a^{C_N} (a^{-1})^{C_N}$ and $e^G= \{e\}$, it follows that $a^{C_N} (a^{-1})^{C_N}=\{e\}$ and so
  $C_N={\bf C}_G(a)$. Since $C_N/N={\bf C}_{\bar{G}}(\bar{a})$, the result follows.
 \end{proof}
 
\begin{proof}[Proof of Theorem B]
We are going to use induction on $\eta(A A^{-1})|G|$. 
Observe that the statement is 
true if $\eta(A A^{-1})=1$, since in that case
 $A A^{-1}=\{e\}$ and thus ${\bf C}_G(A)=G$.

We are going to use the notation of Hypothesis \ref{hypo}, where $N$, in addition, is 
a chief factor of $G$. Thus $A=a^G$ and ${\bf C}_G(A)=\core_G(C)$.
  Lemma \ref{lema2} implies that
 $\eta(\bar{a}^{\bar{G}} ({\bar{a}}^{-1})^{\bar{G}})\leq \eta(a^G (a^{-1})^G)$. Since $|\bar{G}|<|G|$, by induction we have that 
\begin{equation}\label{induction}
 \dl(\bar{G}/\core_{\bar{G}}({\bf C}_{\bar{G}} ( \bar{a})))\leq 2\eta(\bar{a}^{\bar{G}} ({\bar{a}}^{-1})^{\bar{G}})-1.
 \end{equation}
 
Observe that $\bar{G}/\core_{\bar{G}}({\bf C}_{\bar{G}}(\bar{a}))$ is isomorphic to
$G/\core({\bf C}_G(C_N))$. Therefore
\begin{equation}\label{induc}
\dl(G/\core({\bf C}_G(C_N))=\dl(\bar{G}/\core_{\bar{G}}({\bf C}_{\bar{G}}(\bar{a}))).
\end{equation}
  Assume that  $\core_G(C_N)=\core_G(C)$. 
 Then $G/\core_G(C)$ is isomorphic to 
 the group $\bar{G}/\core_{\bar{G}}({\bf C}_{\bar{G}} ( \bar{a}))$.  By \eqref{induction}, \eqref{induc} and 
 Lemma \ref{lema2},
 we have
 \begin{equation*}
 \begin{split}
 \dl(G/\core_G(C)) &=\dl(\bar{G}/\core_{\bar{G}}({\bf C}_{\bar{C}} ( \bar{a}))\\
 &\leq 2\eta(\bar{a}^{\bar{G}} ({\bar{a}}^{-1})^{\bar{G}})-1 \leq 2\eta(a^G (a^{-1})^G)-1.
 \end{split}
 \end{equation*}
 We may assume then that 
  $\core_G(C)$ is a proper subset of $\core_G(C_N)$. 
   Observe that then 
  $C$ is properly contained in $C_N$ and therefore, by Lemma \ref{equal}, we must have
  that $$\eta(\bar{a}^{\bar{G}} ({\bar{a}}^{-1})^{\bar{G}})<\eta(a^G (a^{-1})^G).$$ 
  So $\eta(\bar{a}^{\bar{G}} ({\bar{a}}^{-1})^{\bar{G}})+1 \leq \eta(a^G (a^{-1})^G)$. Since $G$ is 
  supersolvable, $N$ is cyclic.  By Lemma \ref{intersection},
  \eqref{supersolvable}, \eqref{induction} and 
  \eqref{induc} we have that
   \begin{equation*}
 \begin{split}
 \dl(G/\core_G(C))
 & \leq \dl(G/\core_G(C_N)) +  \dl(\core_G(C_N)/\core_G(C))\\
 &= \dl(\bar{G}/\core_{\bar{G}}({\bf C}_{\bar{G}} ( \bar{a})))+  \dl(\core_G(C_N)/\core_G(C)) \\
 &\leq [2\eta(\bar{a}^{\bar{G}} ({\bar{a}}^{-1})^{\bar{G}})-1] +2 \\
 &\leq 2(\eta(a^G (a^{-1})^G)-1)+2 -1+2\\
 &\leq 2 \eta(a^G (a^{-1})^G) -1.
 \end{split}
 \end{equation*}
 The proof is now complete.
 \end{proof}  
\end{section}
\begin{section}{Proof of Corollary C}
\begin{lem}\label{centerintersection}
Let $G$ be a finite group, $A$ and $B$ be conjugacy classes of $G$ such that $AB\cap {\bf Z}(G)\neq \emptyset$. Then  $\eta(AB)=\eta(AA^{-1})$.
\end{lem} 
\begin{proof}
Since $AB\cap {\bf Z}(G)\neq \emptyset$, there exist some
$z\in {\bf Z}(G)$, $a\in A$ and $b\in B$ such that
 $z=ab$. Thus
$b=a^{-1}z$ and so $AB=a^G b^G= a^G (a^{-1}z)^G =(a^G (a^{-1})^G)z=ABz$. It follows then
that  $\eta(AB)=\eta(AA^{-1})$.
\end{proof}

Corollary C follows from Lemma \ref{centerintersection} and Theorem B.
\end{section}

\end{document}